\documentclass[12pt]{article}  
\usepackage{amssymb}
\usepackage{amsmath,amsthm,enumerate}
\theoremstyle{plain}
\newtheorem{thm}{Theorem}[section]
\newtheorem{prop}[thm]{Proposition}
\newtheorem{cor}[thm]{Corollary}
\newtheorem{lemma}[thm]{Lemma}

\theoremstyle{definition}

\newtheorem{rmk}[thm]{Remark}

\newenvironment{enum1}{\begin{enumerate}[\upshape (1)]}{\end{enumerate}}

\title{Zeros of the Hurwitz zeta function\\in the interval $(0,1)$\thanks{Research
    supported by Swiss National Science Foundation Grant no. 107887. Support has also been given over time by Scuola Normale Superiore of Pisa, University of Z\"urich and IBM Z\"urich Research Lab.}}
\author{Davide Schipani\\
Institute of Mathematics\\
University of Zurich\\  
davide.schipani(at)math.uzh.ch}
\date{\today}
\begin{document}\maketitle

\begin{abstract}


We first prove an inequality for the Hurwitz zeta function $\zeta(\sigma,w)$ in the case $\sigma>0$. As a corollary we derive that it has no zeros and is actually negative for $\sigma\in(0,1)$ and $1-\sigma\leq w$  and, as a particular instance, the known result that the classical zeta function has no zeros in $(0,1)$. 

\end{abstract}\medskip

\section{Statement and proof of the results}

The Hurwitz zeta
function is classically defined for $\Re(s)>1$ as 

$$
\zeta(s,w)\doteq\sum_{n=0}^{\infty} (n+w)^{-s},
$$

\noindent with $w$ being a positive real number, and it can be continued analytically to the whole $s$-plane, except for a pole in $1$ (see e.g. \cite[Proposition 9.6.6]{co07}). Sometimes the definition is extended by letting $w$ be a complex number, while in other situations $w$ is only restricted to be a real number in $(0,1]$. Notice in fact the following relation:
$$
\zeta(s,w)=\zeta(s,w+1)+w^{-s},
$$
which follows by considering a summation over the terms $(n+1+w)$.


The following theorem is meant to add a new result to what is already known about its zeros (see e.g. \cite{sp76}). As usual in the literature, $\sigma$ will denote the real part of $s$.

\begin{thm}
Suppose that $\sigma>0$ and $\sigma\neq 1$. Then $\zeta(\sigma,w)<\frac{1-\sigma-w}{(1-\sigma)w^{\sigma}}$.
\end{thm}
\begin{proof}
To make it clear how the steps in the proof drive towards the thesis and where they leave space for possible improvements, we use here the following approach: we start considering $s$ a generic complex number and only later we will add conditions on the parameters.

As a first step we derive a representation for the Hurwitz zeta function through the Euler-Maclaurin summation formula, in analogy with the one derived for the classical Riemann zeta function (see e.g \cite[chapter 6]{ed01}).
Namely, the following Euler-Maclaurin summation formula

$$
\sum_{n=N}^M f(n)=\int_N^M f(x)dx + \frac{1}{2}[f(M)+f(N)]+ \frac{B_2}{2}f'(x)\big|_N^M + \frac{B_4}{4!}f^{(3)}(x)\big|_N^M +\ldots,
$$

\noindent where $B_2, B_4,\ldots$ are Bernoulli numbers, is applied to $f(n)=(n+w)^{-s}$ letting $M$ tend to $\infty$, to obtain, for $\Re(s)>1$,

$$
\zeta(s,w)=\sum_{n=0}^{N-1}(n+w)^{-s}-\frac{(N+w)^{1-s}}{1-s}+O((N+w)^{-\sigma}).
$$


This formula is in fact true for all $\mathbb{C}\backslash\{1\}$, as explained for the analogous case of the Riemann zeta function (\cite[section 6.4]{ed01}) by considering an expression for the remainder in integral form
and its correspondent halfplane of convergence (see also \cite[Proposition 9.6.7]{co07}).

The formula above implies that, if we restrict to $\Re(s)>0$, $s\neq 1$,
$$
\zeta(s,w)=\lim_{N\to\infty}\left(\sum_{n=0}^{N}(n+w)^{-s}-\frac{(N+w)^{1-s}}{1-s}\right).
$$


We notice that

$$\sum_{n=0}^{N}(n+w)^{-\sigma}=w^{-\sigma}+\int_0^N (x+w)^{-\sigma}dx+\lambda_N,$$

\noindent where 



$$\lambda_N=\sum_{n=1}^N\{(n+w)^{-\sigma}-\int_{n-1}^n(x+w)^{-\sigma}dx\}$$

\noindent is such that $0>\lambda_1>\lambda_N>\lambda_{N+1}$.

After computing $\int_0^N (x+w)^{-\sigma}dx= \frac{(N+w)^{1-\sigma}}{1-\sigma}-\frac{w^{1-\sigma}}{1-\sigma}$, we get

$$\sum_{n=0}^{N}(n+w)^{-\sigma}-\frac{(N+w)^{1-\sigma}}{1-\sigma}=\frac{1-\sigma-w}{(1-\sigma)w^{\sigma}}+\lambda_N,$$ 

\noindent so that

$$
\lim_{N\to\infty}\left(\sum_{n=0}^{N}(n+w)^{-\sigma}-\frac{(N+w)^{1-\sigma}}{1-\sigma}\right)=\zeta(\sigma,w)=\frac{1-\sigma-w}{(1-\sigma)w^{\sigma}}+\lim_{N\to\infty} \lambda_N,
$$

\noindent and therefore

$$
\zeta(\sigma,w)<\frac{1-\sigma-w}{(1-\sigma)w^{\sigma}}.
$$






\end{proof}

\begin{cor}
$\zeta(\sigma,w)$ is negative and in particular nonzero when $\sigma\in(0,1)$ and $1-\sigma\leq w$.
\end{cor}


A particular case is the following well known result (see e.g. \cite[chapter 13]{ap76}):

\begin{cor}\label{cor}
The Riemann zeta function $\zeta(s)$ has no zeros in $(0,1)$.
\end{cor}
\begin{proof}
In fact $\zeta(s)=\zeta(s,1)$.
\end{proof}

\begin{rmk}
Notice that when $\sigma>1$ the upper bound is always positive, as it should be, given the definition of $\zeta(\sigma,w)$ for $\Re(s)>1$.
\end{rmk}

\begin{rmk}
It is well known that the Dirichlet $L$-function, defined for $\Re(s)>1$ as
$L(\chi,s)\doteq\sum_{n=1}^{\infty}\frac{\chi(n)}{n^s}$, where $\chi$ is a Dirichlet character (mod $q$), can be represented through a sum of Hurwitz zeta functions (see e.g \cite{co07} or \cite{ap76}) in the following way:
$$
L(\chi,s)=q^{-s}\sum_{a=1}^q \chi(a)\zeta(s,a/q).
$$

The Extended Riemann Hypothesis conjectures that the Dirichlet $L$-function
$L(\chi,s)$, for a
primitive character $\chi$, has no zeros
with real part different from $1/2$ in the critical strip;
and its strong version says that $L(\chi,1/2)$ is always
nonzero too: see also \cite[section 10.5.7]{co07}).

We leave then open for investigation to see whether our main result could help to find new zero-free regions for this class of functions, for example to show that the Dirichlet $L$-functions are also nonzero on the real axis between $0$ and $1$, which would thus 
establish a weaker version of the Extended Riemann Hypothesis, though stronger in including the point $1/2$. A very simple instance is actually at hand:
\begin{cor}
Let $\chi_2$ be the unique character modulo $2$. Then $L(\chi_2,s)$ has no zeros in $(0,1)$.
\end{cor}
\begin{proof}
From the above formula we get that $L(\chi_2,s)=2^{-s}\zeta(s,1/2)$ and we know from our theorem that $\zeta(s,1/2)$ is nonzero in $[1/2,1)$, so that $L(\chi_2,s)$ is also nonzero in $[1/2,1)$. Now, we know by the functional equation for Dirichlet $L$-functions (see e.g. \cite[Theorem 10.2.14]{co07}) that if $s$ is not a zero, then $1-s$ is not a zero either. Therefore $L(\chi_2,s)$ is nonzero throughout $(0,1)$.
\end{proof}

An alternative proof of this result can be derived using Corollary \ref{cor} and the fact that $\zeta(s,1/2)=\zeta(s)\cdot (2^s-1)$ (see e.g. \cite[Proposition 9.6.2]{co07}). In fact this case belongs to a more general scenario: when $\chi$ is a principal character $\chi_0$ modulo $q$, then: 
$$
L(\chi_0,s)=\zeta(s)\prod_{p|q}\left(1-\frac{1}{p^s}\right)
$$ 
(see e.g. \cite{co07} or \cite{sc10}).

\end{rmk}

\section{Acknowledgements}

Thanks are due for suggestions and comments to Joachim Rosenthal and Alessandro Cobbe as well as to the anonymous referee.


\bibliography{huge} \bibliographystyle{plain}

\def\cprime{$'$} \def\polhk#1{\setbox0=\hbox{#1}{\ooalign{\hidewidth
  \lower1.5ex\hbox{`}\hidewidth\crcr\unhbox0}}}
  \def\polhk#1{\setbox0=\hbox{#1}{\ooalign{\hidewidth
  \lower1.5ex\hbox{`}\hidewidth\crcr\unhbox0}}} \def\cprime{$'$}
  \def\cprime{$'$} \def\cprime{$'$} \def\cprime{$'$}
\begin{thebibliography}{1}

\bibitem{ap76}
T.~M. Apostol.
\newblock {\em Introduction to analytic number theory}.
\newblock Springer-Verlag, New York, 1976.

\bibitem{co07}
H.~Cohen.
\newblock {\em Number theory. {V}ol. {II}. {A}nalytic and modern tools}, volume
  240 of {\em Graduate Texts in Mathematics}.
\newblock Springer, New York, 2007.

\bibitem{ed01}
H.~M. Edwards.
\newblock {\em Riemann's zeta function}.
\newblock Dover Publications Inc., Mineola, NY, 2001.
\newblock Reprint of the 1974 original.

\bibitem{sc10}
D.~Schipani.
\newblock Generalized {R}iemann hypotheses: sufficient and equivalent criteria.
\newblock {\em JP Journal of Algebra, Number Theory and Applications},
  19(2):203--214, 2010.

\bibitem{sp76}
R.~Spira.
\newblock Zeros of {H}urwitz zeta functions.
\newblock {\em Math. Comp.}, 30(136):863--866, 1976.

\end{thebibliography}

\end{document}